\documentclass[11pt]{amsart}
\usepackage[french,english]{babel}
\usepackage{amsmath}
\usepackage{amsfonts}
\usepackage{amssymb}
\usepackage{amsthm}
\usepackage{graphicx}
\usepackage[latin1]{inputenc}
\usepackage[T1]{fontenc}
\usepackage{eufrak}
\usepackage[all]{xy}
\newtheorem{thm}{Theorem}[section]
\theoremstyle{definition}
\newtheorem{dfn}{Definition}
\theoremstyle{remark}
\newtheorem{rem}{Remark}
\usepackage{color}
\newenvironment{dem}{\noindent{\textit{Proof.}}}{\begin{flushright}$\diamondsuit$\end{flushright}}
\theoremstyle{plain}
\newtheorem{lem}{\bf{Lemma}}[section]
\theoremstyle{plain}
\newtheorem{prop}{Proposition}[section]

\begin{document}
\author{Guillaume Pouchin}
\def\sym{{\mathfrak{S}}}
\def\Zz{{\mathbb{Z}}}
\def\Nn{{\mathbb{N}}}
\def\Cc{{\mathbb{C}}}
\def\Drap{{\mathfrak{B}}}
\def\Nil{{\mathcal{N}}}
\def\Fq{{\mathbb{F}_{q}}}
\def\Hn{{\mathcal{H}_{n}}}
\def\Hnaffp{{\widehat{\mathcal{H}}_{n}^{+}}}
\def\Hnaff{{\widehat{\mathcal{H}}_{n}}}
\def\Zz{{\mathbb{Z}}}
\def\calO{{\mathcal{O}}}
\def\flagaffr{{\widehat{\mathcal{F}l}_{r}}}
\def\Ww{{\mathcal{W}}}
\def\Ll{{\mathcal{L}}}
\def\Waffn{{\widehat{\mathfrak{S}}_{n}}}
\def\Wn{{\mathfrak{S}_{n}}}
\def\dbar{{\underline{d}}}
\def\Xchap{{\hat{X}}}
\def\Ychap{{\hat{Y}}}
\def\Xchaplus{{(\hat{X} \times \hat{X})^+}}
\def\Ychaplus{{(\hat{Y} \times \hat{Y})^+}}
\def\CalA{{\mathcal{A}}}
\def\CalB{{\mathcal{B}}}
\def\CalC{{\mathcal{C}}}
\def\flagaffn{{\widehat{\mathcal{F}l}_{n}}}
\def\flagaffnd{{\widehat{\mathcal{F}l}_{n,d}}}
\def\Ww{{\mathcal{W}}}
\def\Ll{{\mathcal{L}}}
\def\CalX{{\mathcal{X}}}
\def\CalY{{\mathcal{Y}}}
\def\CalZ{{\mathcal{Z}}}
\def\CalF{{\mathcal{F}}}
\def\CalO{{\mathcal{O}}}
\def\End{{\mathrm{End}}}
\def\Hom{{\mathrm{Hom}}}
\def\Im{{\mathrm{Im}}}
\def\Aut{{\mathrm{Aut}}}
\def\dom{{\mathrm{dom}}}
\def\Ker{{\mathrm{Ker}}}
\def\End{{\mathrm{End}}}
\def\lbar{{\underline{\lambda}}}
\def\Hnaffl{{\widehat{\mathcal{H}}_{n,\lambda}}}
\title[Geometric Schur-Weyl duality]{A geometric Schur-Weyl duality for quotients of affine Hecke algebras}

\begin{abstract}
After establishing a geometric Schur-Weyl duality in a general
  setting, we recall this duality in type A in the finite and affine
  case. We extend the duality in the affine case to positive parts of the
  affine algebras. The positive parts have nice ideals coming from geometry,
  allowing duality for quotients. Some of the quotients of the positive
  affine Hecke algebra are then identified to some cyclotomic Hecke algebras and
  the geometric setting allows the construction of canonical bases.
\end{abstract}
\maketitle
\section{Introduction}

The so-called Schur-Weyl duality is a bicommutant theorem which classically
holds between $GL_d(\Cc)$ and the symmetric group $\Wn$. The first group
acts naturally on $V=\Cc^d$ and diagonally on $V^{\otimes n}$. The symmetric
group acts on $V^{\otimes n}$ by permuting the tensors. The theorem says
that the algebras of these groups are the commutant of each other inside
$\mathrm{End}_\Cc(V^{\otimes n})$ when $d \geq n$. More precisely we have
canonical morphisms:
\[
\phi: \Cc[GL_d(\Cc)] \rightarrow \End_\Cc(V^{\otimes n}), \ \ \ \ \ 
\psi: \Cc[\Wn] \rightarrow \End_\Cc(V^{\otimes n})
\]
such that 
\[
\phi(\Cc[GL_d(\Cc)])=\End_\Wn(V^{\otimes n})
\]
 and
\[
\psi(\Cc[\Wn])=\End_{GL_d(\Cc)}(V^{\otimes n}).
\]

The image of the first map is called the Schur algebra, and the second map
is injective when $d \geq n$.

The case of a base field of arbitrary characteristic has been studied by
numerous autors (\cite {CL}, \cite{Gr}).

 This theorem has many other versions, for affine algebras and for quantum
ones (see \cite{J}).

The aim of this article is to establish geometrically a Schur-Weyl duality
between some quotients of a ``half `` of affine $q$-Schur algebras (which
are themselves quotients of the affine quantum enveloping algebra of
$GL_d$) and some quotients of a ``half'' of the affine Hecke algebra
$\Hnaff$ (of type A), when $ d \geq n$.

  The Schur-Weyl duality between affine $q$-Schur algebras and affine Hecke algebras is already
known for some time and can be expressed nicely by considering convolution
algebras on some flag varieties (see \cite{VV}). We are considering
quotients of subalgebras which arise naturally in this geometric
interpretation: the ideals which we use to define our quotients are just
the functions whose support lies in some closed subvariety.

In fact some of the quotients of Hecke algebras defined this way are
particular cases of cyclotomic Hecke algebras (where all parameters are
equal to zero). One interesting outcome of our construction is the
existence of canonical bases for such algebras. These are simply defined as
the restriction of certain simple perverse sheaves to the open subvarieties
we are considering. The Schur-Weyl duality also provides a strong link with
affine $q$-Schur algebras.

The article is organized as follows: the first part is the geometric
setting needed to have a Schur-Weyl duality. The second part is the
application of the first part in type A for the finite and affine case. In
the third part we show that the geometric Schur-Weyl duality remains when
we restrict ourselves to some subalgebras verifying some conditions. We
apply this in the next part to the affine case by taking positive parts of
our affine algebras. These positive algebras have interesting two-sided ideals
coming from geometry, so the fifth part is dedicated to quotients by such
ideals, in particular we establish the Schur-Weyl duality for our quotient
algebras. The quotients are also identified. The sixth part deals with
the construction of canonical bases of our quotients, using their
construction in terms of intersection complexes.

The last part of the article is the study of the case $d < n$, where only
one half of the bicommutant holds. This answers a question of Green
(\cite{G}) in the affine case.

I deeply thank my supervisor O.Schiffmann for his useful help, comments and his availability. 

\section{Schur-Weyl duality in a general setting}

Let $G$ be a group acting on two sets $X$ and $Y$. Let us assume that we have
the following data:
\begin{itemize}
\item a decomposition $Y=\bigsqcup_{i\in I}Y_i$, where $I$ is a finite set,
\item for each $i \in I$, a surjective $G$-equivariant map,
  $\phi_i : X \rightarrow Y_i$, which has finite fibers of constant
  cardinal $m_i$,
\item An element $\omega \in I$ for which the map $\phi_\omega$ is bijective.
\end{itemize}

We equip the products $X \times Y$, $X \times X$ and $Y \times Y$ with the
diagonal $G$-action.

Let  $A=\Cc_{G} (Y \times Y)$ be the set of $G$-invariant functions which
take non-zero values on a finite number of $G$-orbits, and define $B=\Cc_{G} (X \times
X)$ in the same way. These are equipped with the convolution product
\[
f*g(L,L'') = \sum_{L'} f(L,L')g(L',L'').
\]

The space $C=\Cc_{G} (Y \times X)$ is endowed with
a natural action by convolution of $A$ (resp. $B$) on the left (resp. on the right).

\begin{thm}[Bicommutant Theorem]\label{bicom}

We have:
\[
\mathrm{End}_B(C)=A,
\]
\[
\mathrm{End}_A(C)=B.
\]
\end{thm}

\begin{dem}
 Let's prove the first assertion. Let $P \in
\mathrm{End}_B(C)$.

 From the decomposition $Y=\bigsqcup_{i \in I} Y_i$, we can split $C$
 as a direct sum of vector spaces:  
\[
C= \bigoplus_{i \in I} \Cc_G(Y_i \times X)
\] 

and hence $\End(C)$ as:
\begin{equation}
\label{dec}
\End(C)=\bigoplus_{i,j} \Hom(\Cc_G(Y_i \times X), \Cc_G(Y_j \times X))
\end{equation} 
 The $(i,j)$-component $P'_{(i,j)}=P'$ of $P$ with respect to (\ref{dec}) is the
morphism defined by $P'(f)=1_{\calO_{\Delta(j)}}*P(1_{\calO_{\Delta(i)}}*f)$, where $\Delta(l)$
is the diagonal of $Y_l \times Y_l$.

 The $G$-equivariant surjective map $\phi_i: X \rightarrow Y_i$ gives rise to the
 $G$-equivariant maps:
\[
 Id \times \phi_i: Y_j \times X \rightarrow Y_j
 \times Y_i
\]
\[
 \phi_i \times Id: X \times X \rightarrow Y_i \times X
\]
 From these surjective maps we canonically build the injections:
\[
\psi_i : \Cc_G(Y_j \times Y_i) \hookrightarrow
\Cc_G(Y_j\times X)
\]
\[
\chi_j : \Cc_G(Y_j \times X) \hookrightarrow \Cc_G(X \times X)
\]

 For every $f$ in $\Cc_G(Y_i \times X)$ we have:
\[
P'(f)=P'(m_i^{-1}1_{\Delta(Y_i \times X)} *
\chi_i(f))=m_i^{-1}P'(1_{\Delta(Y_i \times X)}) * \chi_i(f)
\]
where $P'(1_{\Delta(Y_i \times X)}) \in \Cc_G(Y_j \times X)$.
  
  We will now prove that $P'(1_{\Delta(Y_i
\times X)})$ belongs to the image of $\Cc_G(Y_j \times Y_i)$ under $\psi_i$.

 By definition the image of $\psi_i$ in $\Cc_G(Y_j
\times X)$ is the set of functions taking the same values on two orbits of $Y_j
\times X$ which have the same image in $Y_j\times Y_i$. We introduce:
\[
Z_i= \{ (L,L') \in X \times X,
\phi_i(L)=\phi_i(L') \}.
\]

Observe that $P'(1_{\Delta(Y_i\times X)})*1_{Z_i}=P'(1_{\Delta(Y_i\times
X)}*1_{Z_i})=m_iP'(1_{\Delta(Y_i\times X)})$. The following lemma implies
that $P'(1_{\Delta(Y_i \times X)})$ belongs to $\Im(\psi_i)$

\begin{lem}
\[
\Im(\psi_i)= \{ h \in \Cc_G(Y_j \times X), \ h *
1_{Z_i}=m_i h \}
\]
\end{lem}
\begin{dem}
 For the inclusion of the left-hand side in the right-hand side, we can
 write for $h \in \Im(\psi_i)$ and $(L, L') \in Y_j
\times X$ :
\[
h * 1_{Z_i}(L,L')=\sum_{L''} h(L,
L'')1_{Z_i}(L'',L') \\ 
=\sum_{\phi_i(L'')=\phi_i(L')} h(L,L'')= m_i h(L,L')
\]

 For the other inclusion, let $h \in \Cc_G(Y_j \times X)$ be such that
$h*1_{Z_i}=m_i h$. Take $(L,M)$ and $(L,N)$ in $Y_i \times X$ such
that $\phi_j(M)=\phi_j(N)$. Then we have
$1_{Z_i}(L',M)=1_{Z_i}(L',N)$ for every $L'$ in $X$. Then:
\[
m_i h(L,M)=h*1_{Z_i}(L,M)=
\sum_{L'}h(L,L')1_{Z_i}(L',M)= 
\]
\[
\sum_{L'}h(L,L')1_{Z_i}(L',N)=
h*1_{Z_i}(L,N)=m_ih(L,N)
\]
and so $h \in \Im(\psi_i)$.
\end{dem}

Let $g:=\psi^{-1}_i(P'(1_{\Delta(Y_i \times X)})) \in \Cc_G(Y_j \times Y_i)$ 

 So we have for $(L,M) \in Y_j \times X$, $g(L,\phi_i(M))=P'(1_{\Delta(Y_i \times X)})(L,M)$.
 We can now prove that: 
\[
\forall f \in \Cc(Y_i \times X), P'(f)=g*f
\]
 Indeed we have seen that $P'(f)=m_i^{-1}P'(1_{\Delta(Y_i \times X)}) *
\chi_i(f)$. But we have:
\[
m_i^{-1}P'(1_{\Delta(Y_i \times X)}) *
\chi(f)(L,M)=m_i^{-1}\sum_{N\in X}P'(1_{\Delta(Y_i \times
  X)})(L,N)\chi_i(f)(N,M)
\]
\[
=m_i^{-1}\sum_{N\in X}g(L,\phi_i(N))f(\phi_i(N),M)
\]
\[
=\sum_{N'\in Y_i}g(L,N')f(N',M)
\]
\[
=g*f(L,M)
\]
So we have the result for $P'$.

To have it for $P$, it suffices to sum on the orthogonal idempotents.
For every $(i,j) \in I^2$ we have built $g_{(i,j)} \in \Cc_G(Y_j \times Y_i)$
such that $\forall f \in \Cc_G(Y_i \times X)$, $P'_{(i,j)}(f)=g_{(i,j)}
*f$. Let $g= \sum_{i,j} g_{(i,j)}$.

 Then for $f=\oplus_i f_i \in \Cc_G(Y\times X)=\bigoplus_i \Cc_G(Y_i \times
 X)$, we have:
\[
P(f)=\sum_{i,j \in I}
1_{\calO_{\Delta(j)}}*P(1_{\calO_{\Delta(i)}}*f)=\sum_{i,j\in I}
P'_{(i,j)}(f_i)=\sum_{i,j\in I} g_{(i,j)}*f_i \\ =g*f.
\]
\bigskip

We now turn to the second assertion.

Take $P \in \End_A(C)$.

 The projector on $\Cc_G(Y_i \times X)$ parallel to the rest of the sum is
the convolution on the left by the function
$1_{\calO_{\Delta_i}}$, where $\Delta_i$ is the diagonal of
$Y_i \times Y_i$. But $P$ commutes with the action of $A$, so these
subspaces are stable.

 The next lemma will allow us to focus on one such
subspace.
\begin{lem}
The $A$-module $C$ is generated by $\Cc_G(Y_\omega \times X)$.
\end{lem}

\begin{dem}
It is sufficient to verify that for every $f$ in $\Cc_G(Y_i \times X)$, we
have the following formula:
\[
f=\psi_\omega^{-1}(f) * 1_{\Delta(Y_\omega \times X)}
\]
where $\psi_\omega$ is the isomorphism deduced from $\phi_\omega$:
\[
\psi_\omega: \Cc_G(Y_i \times Y_\omega) \rightarrow \Cc_G(Y_i \times X)
\]

\end{dem}

As $C$ is generated as a $A$-module by $\Cc_G(Y_\omega \times
X)$, the endomorphism $P$ is entirely determined by its restriction $P'$
to $\Cc_G(Y_\omega \times X)$. Then we can consider $P' \in
\End_{\Cc_G(Y_\omega \times Y_\omega)}(\Cc_G(Y_\omega \times X))$.

 But the canonical isomorphism $\phi_\omega: Y_\omega \rightarrow X$
allows us to identify $B=\Cc_G(X \times X)$ with $\Cc_G(Y_\omega \times
X)$ and $\Cc_G(Y_\omega \times Y_\omega)$. This way we can see $P'$ as
an element of $\End_B(B)=B$.
\end{dem}

\section{Applications}

\subsection{The linear group}

Let $q$ be a power of a prime number $p$ and $\Fq$ the
finite field with $q$ elements.
We note $G=GL_n(\Fq)$.
In the following everything takes place in a vector space $V$ on
$\mathbb{F}_q$ of dimension $n$. We fix an integer $d \geq n$.

The complete flag manifold $X$ is: 
\[
X=\{ (L_{i})_{1 \leq i
  \leq n} \arrowvert \; L_1 \subseteq L_2 \subseteq \cdots \subseteq L_n=V
  \; , \; \mathrm{dim} L_i =i \}
\]

The partial flag manifold $Y$ of length $d$ is:
\[
Y=\{ (L_{i})_{1 \leq i
  \leq d} \arrowvert \; L_1 \subseteq L_2 \subseteq \cdots \subseteq L_d
  \subseteq V \}
\]

The group $G$ acts canonically on the varieties $X$ and $Y$.

A composition of $n$ of length $d$ is a sequence of integers $\dbar=(d_{1}, \cdots,
d_{d})$ which have a sum equal to $n$. Let $\Lambda(n,d)$ be the set of
compositions of $n$ of length $d$.

For each composition $\dbar$ of $n$ we have a connected component $Y_\dbar$ of
$Y$ defined by $Y_\dbar=\{ (L_\bullet)\in Y, \; \forall i \ \mathrm{dim}
(L_{i+1}/L_i)=d_i \}$ and the decomposition:
\[
Y=\bigsqcup_{\dbar \in \Lambda(n,d)} Y_\dbar 
\]
 Also we have a canonical surjective $G$-equivariant map $\phi_\dbar :
 X \rightarrow Y_\dbar$.

As $d \geq n$, the element $\omega = (\underbrace{1, \cdots
,1}_{n}, 0, \cdots ,0)$ belongs to $\Lambda(n,d)$. The canonical morphism
$\phi_\omega$ is then an isomorphism.

The hypotheses are verified so we can apply the theorem in section 1. In
this case the algebras constructed are well known:

\begin{prop}
The convolution algebra $\Cc_G(X \times X)$ is isomorphic to the Hecke
algebra $\Hn$ with parameter $q=v^{-2}$.
\end{prop}
\begin{prop}
The convolution algebra $\Cc_G(Y \times Y)$ is isomorphic to the $q$-Schur algebra $S_q(n,d)$.
\end{prop}
Thus we obtain the standard Schur-Weyl duality.

\subsection{The affine case}

Let us write $\mathcal{M}$ for the set of $\Fq[[z]]$-submodules of
$(\Fq((z)))^n$ which are free of rank $n$.

The complete affine flag variety $\Xchap$ is:
\[
 \Xchap = \{ (L_{i})_{i  \in \Zz} \in \mathcal{M}^\Zz \arrowvert \  \forall i \:
 L_{i} \subseteq  L_{i+1} , \,  \:  L_{i+n}=z^{-1}L_i , \,
  \: \mathrm{dim}_\Fq L_i/L_{i-1}=1 \}
\]

The affine partial flag variety $\Ychap$ of length $d$ is:
\[
\Ychap = \{ (L_{i})_{i \in \Zz} \in \mathcal{M}^\Zz \arrowvert \  \forall i \:
L_{i} \subseteq L_{i+1} , \, \: L_{i+d}=z^{-1}L_i \}
\]

 Let $G$ be $GL_n(\Fq((z)))$. The varieties defined above are equipped with
 a canonical action of $G$.

We still have a decomposition of $\Ychap$:
\[
\Ychap= \bigsqcup_{\dbar \in \Lambda(n,d)} \Ychap_\dbar
\]
where $\Ychap_\dbar$ is the subvariety of $Y$ defined by:
\[
\Ychap_\dbar= \{ (L_\bullet) \in \Ychap, \  \forall i \in \{1,\cdots,d \} \
\mathrm{dim}_\Fq(L_i/L_{i-1})=d_i \}
\]

 For each element $\dbar \in \Lambda(n,d)$, we have a $G$-equivariant
 surjective map:
\[
\phi_\dbar:\Xchap \rightarrow \Ychap_\dbar
\]
defined by $\phi_\dbar(L_\bullet)_0=L_0$ and
$\phi_\dbar(L_\bullet)_i=L_{\sum_{k=1}^i d_k}$ for $i
\in \{1, \cdots,d\}$.

We can identify these algebras as we did in the previous paragraph:
\begin{prop}[\cite{IM}]
The algebra $\Cc_G(\Xchap \times \Xchap)$ is isomorphic to the affine Hecke
algebra $\Hnaff$ with parameter $q=v^{-2}$.
\end{prop}

\begin{prop}[\cite{VV}]
The algebra $\Cc_G(\Ychap \times \Ychap)$ is isomorphic to the affine $q$-Schur
algebra $\hat{S}_q(n,d)$.
\end{prop}

\section{Subalgebras deduced from subvarieties}

We will now see that the bicommutant theorem remains true, under some
additionnal hypothesis, for a subspace of $\Cc_G(Y\times X)$ and
subalgebras of $\Cc_G(Y\times Y)$ and $\Cc_G(X \times X)$.

Let $X,Y,Y_i$ be as in section 1.

Suppose we are given $G$-subvarieties $\CalZ \subseteq Y \times X$, $\CalX
\subseteq X \times X$, $\CalY \subseteq Y \times Y$ satisfying the
following conditions:
\begin{itemize}
\item for every $i,
j \in I$, when we write $\CalZ_i=\CalZ \cap (Y_i \times X)$ and
$\CalY_{i,j}=\CalY \cap (Y_i \times Y_j)$, then
\[
(\phi_i \times Id_X) (\CalX) = \CalZ_i
\]
and
\[
(\phi_i \times \phi_j) (\CalX) = \CalY_{(i,j)}
\]
\item $\Delta_X \subseteq \CalX$ and $\Cc_G(\CalX)$ is a subalgebra of
  $\Cc_G(X \times X)$.
\end{itemize}

From the above assumptions it follows that:

\begin{itemize}
\item $\Cc_G(\CalZ)$ is stable for the action of $\Cc_G(\CalX)$ on $\Cc_G(Y \times X)$.
\item $\Cc_G(\CalY)$ is a subalgebra of $\Cc_G(Y \times Y)$, and
  $\Cc_G(\CalZ)$ is stable for the action of $\Cc_G(Y \times X)$.
\item the spaces $\Cc_G(\CalX)$, $\Cc_G(\CalY)$ and $\Cc_G(\CalZ)$
  contain the characteric functions of the diagonals $1_{\CalO_{\Delta(X \times X)}}$,
  $1_{\CalO_{\Delta(Y_i \times X)}}$ and $1_{\CalO_{\Delta(Y_i \times Y_j)}}$
  (for every $i$, $j$ in $I$).
\item the subspace $\Cc_G(\CalZ_\omega)$ generates $\Cc_G(\CalZ)$ as a
  $\Cc_G(\CalY)$-module.
\item the diagonal function $1_{\CalO_{\Delta(Y_i \times X)}}$ generates
  $\Cc_G(\CalZ_i)$ as a $\Cc_G(\CalX)$-module.
\item we have isomorphisms deduced from $\psi$ and $\chi$:
\[
\Cc_G(\CalY_{(\omega,\omega)}) \simeq \Cc_G(\CalZ_\omega) \simeq \Cc_G(\CalX)
\]
\end{itemize}

\begin{thm}\label{th3}
Under the previous conditions, the following bicommutant theorem holds:
\[
\End_{\Cc_G(\CalX)}(\Cc_G(\CalZ))=\Cc_G(\CalY)
\]
\[
\End_{\Cc_G(\CalY)}(\Cc_G(\CalZ))=\Cc_G(\CalX)
\]
\end{thm}

The proofs are the same as in the case of the whole space.

\section{The positive part of the affine Hecke algebra}

We get back to the setting of section 2.2. Thus $\Xchap$ and $\Ychap$ are
resp. the complete affine flag variety and the partial affine flag
variety. We recall that we take $d \geq n$, where $n$ is the rank of the
free modules and $d$ is the periodicity in the partial affine flag
variety. Consider the subvarieties:
\begin{itemize}
\item $\CalX= (\Xchap \times\Xchap)^+ = \{ (L_\bullet, L_\bullet') \in \Xchap \times
  \Xchap, L_0' \subseteq L_0\}$
\item $\CalY= (\Ychap \times\Ychap)^+ = \{ (L_\bullet, L_\bullet') \in \Ychap \times
  \Ychap, L_0' \subseteq L_0\}$
\item $\CalZ= (\Ychap \times\Xchap)^+ = \{ (L_\bullet, L_\bullet') \in \Ychap \times
  \Xchap, L_0' \subseteq L_0\}$
\end{itemize}
which give rise to the convolution algebras
\begin{itemize}
\item $A_+=\Cc_G((\Ychap\times \Ychap)^+)$
\item $B_+=\Cc_G((\Xchap\times \Xchap)^+)$.
\end{itemize}
The subspace
\[
C_+=\Cc_G((\hat{Y}\times \hat{X})^+)
\]

 is a $(A_+,B_+)$-bimodule. It is easy to check that the hypothesis of
 theorem 3.1 are verified, so that the bicommutant theorem still holds:
\begin{prop}
We have:
\[
\End_{A_+}(C_+)=B_+
\]
\[
\End_{B_+}(C_+)=A_+
\]
\end{prop}

\begin{dem}
It is a direct application of theorem \ref{th3}.

\end{dem}

Our immediate aim is to identify precisely the algebra $B_+$. For this, we
need to recall in more details the structure of affine Weyl group in type A.

\subsection{The extended affine Weyl group in type A}
Let us first recall first the definition of the extended affine Weyl group  in the
general case of a connected reductive group $G$ over $\Cc$. We write $T$
for a maximal torus of $G$, $W_0 = N_G(T)/T$ is the Weyl group of $G$. The
group $W_0$ acts on the character group $X=\Hom(T,\Cc^*)$, which allows us
to consider the semi-direct product $W=W_0 \rtimes
X$, which is called the extended affine Weyl group of $G$. The root system
$R$ of $G$ generates a sublattice of $X$, noted $Y$. The semi-direct
product $W'=W_0 \rtimes Y$ is called the affine Weyl group of $G$. It's a
Coxeter group, unlike the extended affine Weyl group. It is also a normal
subgroup of $W$.

There is an abelian subgroup $\Omega$ of $W$ such that
$\omega^{-1}S\omega=S$ for every $\omega \in \Omega$ and $W=\Omega \rtimes
W'$.

In the case of $G=GL_n(\Cc)$, the Weyl group $W_0$ is isomorphic to the
symmetric group $\Wn$. We write $S=\{s_1, \cdots,s_{n-1}\}$ for its simple
reflections. The group $W'$ is still a Coxeter group, which is generated by
the simple reflections of $W_0$ and an additional elementary reflection
$s_0$. The group $\Omega$ is isomorphic to $\Zz$, and is generated by an
element $\rho$ which verifies $\rho^{-1}s_i\rho=s_{i+1}$ for every
$i=1,\cdots,n-1$, where we write $s_n$ for $s_0$.

 The group $W'$ is then the group generated by the elements
$s_0,\cdots,s_{n-1}$, with the following relations:
\begin{enumerate}
\item $s_i^2=1$ for every $i=1,\cdots n-1$
\item $s_i s_{i+1} s_i = s_{i+1} s_i s_{i+1}$ where the indices are taken
  modulo $n$.

\end{enumerate}

We have $W=\Omega \rtimes W'$, where the group $\Omega$ is isomorphic to
$\Zz$, generated by an element $\rho$ which verifies:
\[
\rho^{-1}s_i\rho=s_{i-1}.
\]

The character group of a torus of $GL_n(\Cc)$ is naturally
isomorphic to $\Zz^n$. Thus we have $W=\Wn \ltimes \Zz^n$, where the group
$\Wn$ acts on $\Zz^n$ by permutation.

The group $W=\Wn \ltimes \Zz^n$ can also be considered as a subgroup of the
group of the automorphisms of $\Zz$ in the following way: to each $(\sigma,
(\lambda_i)) \in \Wn \ltimes \Zz^n$, we associate the element
$\tilde{\sigma} \in \Aut(\Zz)$ defined by: 
\[
\tilde{\sigma}(i)=\sigma(r) + kn + \lambda_{r}n
\]
where $i=kn+r$ is the Euclidian division of $i$ by $n$, taking the rest
between $1$ and $n$.

In fact if we write $\tau$ for the element of $\Aut(\Zz)$ defined by:
\[
\tau : i \mapsto i + n
\]
and set $\Aut_n(\Zz)= \{ \sigma \in \Aut(\Zz), \  \sigma \tau = \tau \sigma
\}$. Then we obtain the following isomorphism:
\begin{lem}
The previous map provides an isomorphism of groups:
\[
\Wn \ltimes \Zz^n \simeq \Aut_n(\Zz)
\]
\end{lem}

Under this isomorphism, the element $s_i$ ($0 \leq i \leq n-1$) is mapped to $\tilde{s}_i$ defined by:
\[
\left\{
\begin{array}{ll}

\tilde{s}_i(j) = j  & \qquad \mathrm{if} \quad  j \neq i,i+1  \  mod(n), \\
      \tilde{s}_i(j)=j+1 & \qquad \mathrm{if} \quad j= i \  mod(n), \\
       \tilde{s}_i(j)=j-1 & \qquad \mathrm{if}  \quad j =i+1 \  mod(n) .

\end{array}
\right.
\]

The element $\rho$ is mapped to $\tilde{\rho}$ defined by $\tilde{\rho}(i)=i+1$.

The orbits of the action of $G$ on $\Xchap \times \Xchap$
are parametrised by the elements of the extended affine Weyl group
$\Waffn$. Then we can write $\calO_w$ for an orbit, with $w$ in
$\Waffn$. This can be done explicitly in the following way. A couple of
flags $L_\bullet$ and $L'_\bullet$ are in the orbit $w$ if there is a
base $e_1, \cdots,e_n$ of the $\Fq((z))$-module $\Fq((z))^n$ such that 
\[
L_i= \prod_{w(j)\leq i} \Fq e_j \ \ \ \text{and} \ \ \  L'_i=\prod_{j\leq i} \Fq e_j 
\]
where we define $e_i$ for all $i \in \Zz$ by the condition $e_{i+kn}=z^{-k}e_i$
for all $k\in \Zz$.

\begin{thm}[\cite{IM}]
The algebra $\Cc_{G} (\Xchap \times \Xchap)$ is isomorphic to the affine Hecke
algebra $\Hnaff$ specialized at $v^{-2}=q$ and the isomorphism is given by:
\[
\phi: 1_{\calO_{w}} \mapsto T_{w}
\]
for every $w \in \Waffn$.
\end{thm}

To identify the positive part of the affine Hecke algebra, it is necessary to
recall its different presentations.

\subsection{The affine Hecke algebra}\label{aff}

\begin{dfn}[The affine Hecke algebra $\Hnaff$]
The affine Hecke algebra $\Hnaff$ is a $\Cc[v,v^{-1}]$-algebra which may be
defined by generators and relations in either of the following ways:
\\
(1) The generators are the $T_w$, for $w \in \Waffn=\mathfrak{S}_n \ltimes
\Zz^n$.
    The relations are:
\begin{enumerate}
\item $T_{w}T_{w'}=T_{ww'}$     if    $l(ww')=l(w)+l(w')$,
\item $(T_{s_i}+1)(T_{s_i}-v^{-2})=0$ for $s_i=(i,i+1)$.
\end{enumerate}
(2) The generators are $T_i^{\pm 1}$, $i=1 \ldots n-1$ and $X_j^{\pm 1}$,
$j=1 \ldots n$.
    The relations are:
\begin{enumerate}
\item $T_{i}T_{j}=T_{j}T_{i}$    if   $|i-j|>1$,
\item $T_{i}T_{i+1} T_{i}=T_{i+1} T_i T_{i+1}$,
\item $T_i T_i^{-1}=T_i^{-1} T_i=1$,
\item $(T_i+1)(T_i-v^{-2})=0$,
\item $X_i X_i^{-1}=X_i^{-1} X_i=1$,
\item $X_i T_j=T_j X_i$ if $i\neq j,j+1$,
\item $T_i X_i T_i=v^{-2} X_{i+1}$.
\end{enumerate}
\end{dfn}

The isomorphism $\psi$ between these two presentations is uniquely defined
by the following conditions:
\[
   \psi(T_{s_i})= T_{i}
\]
\[
    \psi(\tilde{T}^{-1}_{(\lambda_1,\cdots,\lambda_n)})=
    X_1^{\lambda_1}\cdots X_n^{\lambda_n}
\]
if $\lambda$ is dominant, which means $\lambda_1 \geq \lambda_2 \geq
\cdots \geq \lambda_n$, and if we write $\tilde{T}_w=v^{-l(w)}T_w$, where $l(w)$
is the length of $w$.

One checks that:
\[
T_\rho \mapsto v^{1-n} X_1^{-1}T_1\cdots T_{n-1}.
\]

The multiplication map defines an isomorphism of $\Cc$-vector spaces
\[
\Hnaff \simeq \Cc[\Wn] \otimes_{\Cc} \Cc[v,v^{-1}][X_1^{\pm 1},\ldots,X_n^{\pm 1}]
\]

\subsection{The positive subalgebra}
We are now ready to describe the convolution algebra $B^+$ of $G$-invariant
functions on the positive part $\CalX$ of the product variety $\Xchap
\times \Xchap$.

\begin{thm}
The algebra $\Cc_{G}(\CalX)$ is isomorphic to the subalgebra $\Hnaffp$ of $\Hnaff$
generated by $\Hn$ and the elements $X_{i}$. This means that as a vector space, we have:
\[
\Cc_{G}(\CalX) \simeq \Cc[\Wn] \otimes_{\Cc} \Cc[v,v^{-1}][X_1, \cdots , X_{n}]
\]
\end{thm}

\begin{dem} 

The first observation is that every element $T_i$ is in $\Cc_G(\CalX)$.

The element $X_{1}$ is in $\Cc_{G}(\CalX)$ because, by the isomorphism $\psi$
introduced section \ref{aff}, we have $X_1= v^{1-n} T_{n-1} \cdots T_1
T_\rho^{-1}$. But the element $T_\rho^{-1}$ is the characteristic function
of an orbit in $\CalX$. As $X_{1}$ is in
$\Cc_{G}(\CalX)$, the relations (7) prove that the $X_{i}$s
are in $\Cc_{G}(\CalX)$ as well.

We will now prove that the algebra $\Cc_G(\CalX)$ is generated by the elements
$T_\rho^{-1}$, $ T_1,\cdots , T_{n-1}$. For this purpose we first get
back to the groups.

\begin{lem}
The subsemigroup $\Wn \ltimes \Zz_-^n$ of $\Wn \ltimes \Zz^n$ is generated
by the elements $\rho^{-1}, s_1, \cdots ,s_{n-1}$.
\end{lem}

\begin{dem}
First it is clear that the elements $s_1, \cdots,s_n$ are in $\Wn
\ltimes \Zz_-^n$, because they are in $\Wn$. The element $\rho^{-1}$ can be
written in $\Wn \ltimes \Zz^n$ as $((n\cdots 21),(0,\cdots ,0,-1))$ (to see
it, use the bijection $\Wn \ltimes \Zz^n \simeq \Aut_n(\Zz)$). This element
belongs to $\Wn \ltimes \Zz_-^n$.

We now prove that every element $w$ of $\Wn \ltimes \Zz_-^n$ can be
written as a product of elements among $\rho^{-1}, s_1,\cdots,s_{n-1}$.

We define the degree of an element $w=(\sigma, (\lambda_i))\in \Wn
\ltimes \Zz^n$ by $d=\sum_{i=1}^n \lambda_i$. Let's prove the result by
induction on the degree $d$ of $w\in \Wn \ltimes \Zz_-^n$.

If $d=0$, the result is true because $w$ is an element in $\Wn$.

For $d<0$, let's consider $w'=\rho w$. The degree of $w'$ is $d+1$, so we
have the result by recurrence if $w' \in \Wn \ltimes \Zz_-^n$. Only the
case where $w' \notin \Wn \ltimes \Zz_-^n$ remains.

Under the isomorphism $\Wn \ltimes \Zz^n\simeq \Aut_n(\Zz)$, the subset
$\Wn \ltimes \Zz_-^n$ is mapped to $\{ \tilde{s} \in \Aut_n(\Zz),\  \forall
i=1, \cdots, n,\  \tilde{s}(i) \leq n \}$. The conditions $w\in \Wn \ltimes
\Zz_-^n$ and $w'=\rho w \notin \Wn \ltimes \Zz_-^n$ give in
$\Aut_n(\Zz)$: for every $i=1,\cdots,n$, $\tilde{w}(i)\leq n$ and $\exists j$,
$1\leq j \leq n$, $\tilde{w}'(j)=\tilde{w}(j)+1 \geq n+1$. For this $j$
we have: $\tilde{w}(j)+1=n+1$ hence $\tilde{w}(j)=n$, which is equivalent to
$\sigma(j)=n$ and $\lambda_j=0$.

Besides, we have $d<0$, so there is a $k$ such that $\lambda_k <0$. We
write $t\in \Wn$ for the transposition $(kj)$ and we consider $w''=wt$. Then
$w''=(\sigma t, (\lambda_i)_1^n)$ and so for every $i=1,\cdots ,n$,
$\tilde{w}''(i)<n$. We just saw that this is equivalent to $\rho w'' \in
\Wn \rtimes \Zz_-^n$. So we can apply our recurrence to $x=\rho w''$, which
is of degree $d+1$, to obtain that $x$ is in the subsemigroup generated by
$\rho^{-1}, s_1, \cdots ,s_{n-1}$. As $w=\rho^{-1} xt$, this is also true for $w$.
\end{dem} 

To lift this result to the Hecke algebra, we need a little more: we have to prove
that every element of $\Wn \ltimes \Zz_-^n$ has a reduced decomposition as a
product of $s_1,\cdots,s_{n-1},\rho^{-1}$.

\begin{lem}\label{redgroup}
Every element of $\Wn \ltimes \Zz_-^n$ has a reduced decomposition which
involves only the elements $s_1,\cdots,s_{n-1},\rho^{-1}$.
\end{lem}

\begin{dem}
For $w\in \Wn \ltimes \Zz_-^n$, we write $k$ for its length and $d$ for its
degree ($d \leq 0$). We now proceed by induction on $k-d$.

If $k=0$, the element $w$ is of length $0$ so it is a power of $\rho$,
which is negative because $w\in\Wn \ltimes \Zz_-^n$. We are done.

Si $d=0$, the element $w$ is of degree $0$ in $\Wn \ltimes \Zz_-^n$ so it
is an element of $\Wn$. Then we are done  because an element of $\Wn$ has
a minimal decomposition which uses only $s_1,\cdots,s_{n-1}$.

The last case is when $d<0$ and $k>0$. We know that $w$ has a minimal
decomposition of the form:
\[
w=\rho^l s_{i_1}\cdots s_{i_k}
\]
where $0 \leq i_r \leq n-1$.

We now split the proof in two cases:

If $s_{i_k}\neq s_0$, then we can apply the recurrence to $ws_{i_k}$, which
has the same degree as $w$ and whose length is $l(w)-1$. We deduce from
this a minimal writing of $ws_{i_k}$ which involves only the elements
$s_1,\cdots s_{n-1},\rho^{-1}$, then a minimal writing of $w$
using only these elements, by multiplying by $s_{i_k}$.

The remaining case is when $s_{i_k}=s_0$. In this case we have
$l(ws_0)<l(w)$. We know at this point by using the isomorphism $\Wn
\ltimes \Zz_-^n \simeq \Aut_n(\Zz)$(cf \cite{S}, Cor 4.2.3 or \cite{G}, Cor
1.3.3)) that
$l(ws_0)<l(w)$ implies that $\tilde{w}(0)>\tilde{w}(1)$. Take the
element $w'=ws_0\rho$, and let's show that it belongs to $\Wn \ltimes
\Zz_-^n$. We need to show that for every $i$ such that $1\leq i \leq
n$, we have $\tilde{w}'(i) \leq n$.

By definition $w'(i)=ws_0(i+1)$. So if $1 \leq i \leq  n-2$, we have that
$\tilde{w}'(i)=\tilde{w}(i+1) \leq n$ because $w \in \Wn \ltimes
\Zz_-^n$. We have $\tilde{w}'(n)=\tilde{w}(s_0(n+1))=\tilde{w}(n) \leq n$
too because $w \in \Wn \ltimes \Zz_-^n$. We can deduce that
$\tilde{w}'(n-1)=\tilde{w}(n+1)=n+\tilde{w}(1)<n+\tilde{w}(0)=\tilde{w}(n)
\leq n$, from which it follows that $w' \in \Wn \ltimes \Zz_-^n$. 

As $w'$ has a length equal to $k-1$ and degree $d+1$ and is in the
semigroup $\Wn \ltimes \Zz_-^n$, we can apply the recurrence hypothesis:
$w'$ has a minimal writing which involves only
$s_1,\cdots,s_{n-1},\rho^{-1}$. But as $w=w'\rho^{-1}s_0=w's_1\rho^{-1}$,
we have a minimal writing of $w$ using the $s_1,\cdots,s_{n-1},\rho^{-1}$.
\end{dem}

Observe that by construction, we have:
\[
\Cc_G(\CalX)=\bigoplus_{w \in \Wn\ltimes \Zz_-^n}
\Cc[v,v^{-1}]1_{\CalO_w}=\bigoplus_{w \in \Wn\ltimes \Zz_-^n}\Cc[v,v^{-1}]T_w
\]
By lemma \ref{redgroup}, any $T_w$ may be written as a product of elements \linebreak
$T_\rho^{-1},T_1,\cdots,T_{n-1}$. We easily check that the algebra
generated by these elements is precisely $\Cc[\Wn] \otimes_\Cc \Cc[v,v^{-1}][X_1,\cdots,X_n]$.

\end{dem}

\section{Quotients}

Now that we have at our disposal a bicommutant theorem for the positive
parts of the Hecke algebra and the Schur algebra, we can try to find
subvarieties whose corresponding subalgebras are two-sided ideals of these
algebras, which allows us to take quotients and hope to still have a
bicommutant theorem. 

Let $\lambda=(\lambda_i)_{i=1}^n \in \Nn^n$ be a dominant partition
(i.e. $\lambda_1 \geq \cdots \geq \lambda_n$). For every $(L_\bullet,L_\bullet')\in
\CalX$, as $L_0' \subseteq L_0$ are two free $\Fq[[z]]$-modules of rank
$n$, the quotient $L_0/L_0'$ is a torsion $\Fq[[z]]$-module of rank at most
$n$. We know that the isomorphism classes of torsion $\Fq[[z]]$-modules of
rank at most $n$ are parametrized by the dominant $n$-weights
$\mu_1,\cdots,\mu_n$, with $\mu_1 \geq \cdots \geq \mu_n$.

We define:
\[
\CalX_\lambda= \{ (L_\bullet,L_\bullet') \in \CalX , \
L_0/L_0' \text{ of type }\mu
\text{, with } \forall i=1, \cdots, d,\  \lambda_i \leq \mu_i \}
\]

\begin{lem}
The set $\Cc_G(\CalX_\lambda)$ is a two-sided ideal of $\Cc_G(\CalX)$.
\end{lem}
\begin{dem}
We have to check that if $f$ and $g$ are in $\Cc_G(\CalX)$ with $f$
supported on $\CalX_\lambda$, then $f*g$ and $g*f$ are supported on
$\CalX_\lambda$.

But if $(L_\bullet,L_\bullet') \in \CalX$ and $(L_\bullet',L_\bullet'') \in \CalX_\lambda$, we have
that $L_0'' \subseteq L_0' \subseteq L_0$ and $L_0/L_0'$ is of type $\mu$,
with $\forall i= 1, \cdots, n$, $\lambda_i\leq \mu_i$. From the inclusion
$L_0'/L_0'' \subseteq L_0/L_0''$ we deduce that $L_0/L_0''$ is of type
$\nu$ with $\nu_i \geq \mu_i $ $\forall i$. So $\nu_i \geq \lambda_i$ and
$(L_.,L_.'')$ belongs to $\CalX_\lambda$. Finally $f*g$ is supported on
$\CalX_\lambda$.

If $(L_\bullet,L_\bullet') \in \CalX_\lambda$ and $(L_\bullet',L_\bullet'') \in
\CalX$ then we have $L_0'' \subseteq L_0' \subseteq L_0$ and
$L_0/L_0''$ is of type $\mu$ with $\mu_i\geq \lambda_i$ $\forall i=1 \cdots
n$. As $L_0/L_0'$ is a quotient of $L_0/L_0''$, the type $\nu$ of
$L_0/L_0''$ verifies $\nu_i \geq \mu_i$. So $\nu_i \geq \lambda_i$ and
$(L_\bullet,L_\bullet'') \in \CalX_\lambda$, which gives that $g*f$ is supported on
$\CalX_\lambda$.

\end{dem}

We now define for $i,j$ in $I$:
\[
\CalZ_{\lambda,i} = (\phi_i \times Id) (\CalX_\lambda) \subseteq \CalZ_i
\]
\[
\CalY_{\lambda,i,j} = (\phi_i \times \phi_j) (\CalX_\lambda) \subseteq \CalY_{i,j}
\]
\[
\CalZ_\lambda= \bigsqcup_{i \in I} \CalZ_{\lambda,i}
\]
\[
\CalY_\lambda= \bigsqcup_{i,j \in I} \CalY_{\lambda,i,j}
\]

In the same way $\Cc_G(\CalX_\lambda)$ is a two-sided ideal of
$\Cc_G(\CalX)$, we have the following statement:
\begin{lem}
The set $\Cc_G(\CalY_\lambda)$ is a two-sided ideal of $\Cc_G(\CalY)$.
\end{lem}

The actions of $\Cc_G(\CalX_\lambda)$ and $\Cc_G(\CalY_\lambda)$ map the
space $\Cc_G(\CalZ)$ to $\Cc_G(\CalZ_\lambda)$. Now write
$A_\lambda=\Cc_G(\CalY)/\Cc_G(\CalY_\lambda)$, 
$B_\lambda=\Cc_G(\CalX)/\Cc_G(\CalX_\lambda)$ and
$C_\lambda=\Cc_G(\CalZ)/\Cc_G(\CalZ_\lambda)$. The quotient space
$C_\lambda$ is then a $(A_\lambda,B_\lambda)$-bimodule. We will write
$\Hnaffl$ for $B_\lambda$ in the next part.

We can now state:

\begin{thm}[Bicommutant of the quotient]
\[
\End_{A_\lambda}(C_\lambda)=B_\lambda
\]
\[
\End_{B_\lambda}(C_\lambda)=A_\lambda
\]
\end{thm}
\begin{dem}
We prove the first assertion.

Let $P \in \End_{A_\lambda}(C_\lambda)$. 
As in the theorem \ref{bicom}, the fact that the endomorphism $P$ commutes with
the action of $A_\lambda$ implies that it commutes with the action of the
projectors on $C_{\lambda,i}$ where
$C_{\lambda,i}=\Cc[\CalZ_i]/\Cc_{G}[\CalZ_{\lambda,i}]$, so the subspaces
$C_{\lambda,i}$ are stable by $P$.

We also know that as an $A_\lambda$-module, $C_\lambda$ is generated by
$C_{\lambda,\omega}$. So we only have to study the restriction of $P$ to
$C_{\lambda,\omega}$, where $P$ commutes with the action of
$A_{\lambda,\omega,\omega}$.

But there are canonical isomorphisms $A_{\lambda,\omega,\omega} \simeq
C_{\lambda,\omega} \simeq B_\lambda$. Then we can consider that $P$ belongs
to $\End_{B_\lambda}(B_\lambda)$ which is equal to $B_\lambda$. This proves the
result.

Now we prove the second point of the theorem.

Let $P$ be in $\End_{B_\lambda}(C_\lambda)$, and let $P_i$ be its restriction to
$C_{\lambda,i}$, so that $P_i \in \Hom_{B_\lambda}(C_{\lambda,i},C_\lambda)$.

The canonical morphisms given in theorem \ref{bicom} go through to the positive
parts to give for every $i \in I$ an injection:
\[
\chi_i: C_i^+ \hookrightarrow B^+
\]
whose left inverse is given by the left multiplication by
$\frac{1}{m_i}1_{\Delta_i}$, where $\Delta_i$ is the diagonal in $(\Ychap_i
\times \Xchap)^+$.

We write $\alpha$ for the map from $C_\lambda$ to $C_+$ which associates to
a function in $C_\lambda$ its unique representative in $\Cc_G(\CalZ)$
whose restriction to $\CalZ_\lambda$ is zero.

Now define $P_i' \in \Hom_{B_+}(C_{+,i},C_+)$ by the formula:
\[
P_i'(f)=\alpha(P_i(\frac{1}{m_i}1_{\Delta_i}))*\chi_i(f)
\]
And $P'=\bigoplus_{i\in I} P_i$.

We easily check that $P'$ commutes with the action of $B_\lambda$ and
secondly that through the canonical map $\End_{B_+}(C_+)
\hookrightarrow \End_{B_\lambda}(C_\lambda)$ the morphism $P'$ maps to $P$.

We have lifted $P$ and got $P' \in \End_{B_+}(C_+)$. The bicommutant
theorem for the positive parts gives us the fact that $P'\in A_+$, hence $P
\in A_\lambda$ as desired. We are done.
\end{dem}

We can now identify the two-sided ideal in question.

\begin{prop}
The two-sided ideal $\Cc_G(\CalX_\lambda)$ is generated by the element
$X^{\lambda'}=\prod_{i=1}^{n}X_i^{\lambda_{n-i}}$.
\end{prop}

\begin{dem}
From the definition of $\Cc_G(\CalX_\lambda)$, it is obvious that as a
vector space we can write:
\[
\Cc_G(\CalX_\lambda)=\bigoplus_{\substack{\sigma \in \Wn \\
    \dom(\mu)\geq \lambda}} \Cc T_{(\sigma,-\mu)}
\]
where $\dom(\mu)$ is the partition deduced from $\mu$ by
reordering, and where the partial order between compositions is given by
$\lambda\geq\mu \Leftrightarrow \forall i=1, \cdots, n, \ \lambda_i\geq\mu_i$. In
particular the element $\prod_{i=1}^n X_i^{\lambda_{n-i}}=v^{-l(\lambda)}T_{(Id,-\lambda')}$
belongs to $\Cc_G(\CalX_\lambda)$, and thus $\Hnaff X^{\lambda'}\Hnaff
\subseteq \Cc_G(\CalX_\lambda)$.

The inclusion of the right-hand side in the left-hand side is done.

For the other inclusion, prove first:
\begin{lem} \label{lem1}
For every dominant composition $\nu$ the following holds:
\[
\Hn T_{(Id,\nu)} \Hn = \bigoplus_{\sigma,\sigma' \in \Wn} \Cc T_{(\sigma,\nu^{\sigma'})}
\]
\end{lem}

\begin{dem}
As the element $T_{(Id,\nu)}$ belongs to the sum on the right-hand side
and this space is stable by the action of $\Hn$ on the right and on the
left, the inclusion of the left-hand side in the right-hand side is clear.

We show by induction on the length of $\sigma'$ that every element in the
right-hand side belong to $\Hn T_{(Id,\nu)} \Hn$.

For $\sigma'=Id$: as we have $l(\sigma,
\nu)=l(\sigma,0)+l(Id,\nu)$ because $\nu$ is dominant,
the equation $T_{(Id,\nu)}.T_{(\sigma,0)}=T_{(\sigma,\nu)}$ holds, which
implies that $T_{(\sigma,\nu)} \in \Hn T_{(Id,\nu)} \Hn$.

For $l(\sigma') >0$: we write $\sigma'=ts$, where $s \in S$ and
$l(t)=l(\sigma')-1$.
 By induction we know that for every $u \in \Wn$, we have $T_{(u,\nu^t)}
 \in \Hn T_{(Id,\sigma)} \Hn$. The following holds:
\[
T_{(u,\nu^t)}.T_{(s,0)}= \left\{
                     \begin{array}{ll}
                     T_{(us,\nu^{\sigma'})} & \qquad \mathrm{if} \quad
                     l(us,\nu^{\sigma'})=l(u,\nu^t)+1 \\
                     (1-q)T_{(u,\nu^t)} + qT_{(us,\nu^{\sigma'})} & \qquad
                     \mathrm{if} \quad
                     l(us,\nu^{\sigma'})=l(u,\nu^t)-1 \\
                     \end{array}
                     \right.
\]
As the left-hand side term and $T_{(u,\nu^t)}$ belong to
$\Hn T_{(Id,\nu)} \Hn$ for every $u\in \Wn$, we have also
$T_{(\sigma,\nu^{\sigma'})} \in \Hn T_{(Id, \nu)} \Hn$ for every
$\sigma \in \Wn$.

\end{dem}

To prove the proposition, we first remark that we have the equality \linebreak
$T_{(Id,-\dom(\mu)')}=T_{(Id,-\lambda')}.T_{(Id,-\dom(\mu)'+\lambda')}$ for each
composition $\mu$ such that \linebreak $\dom(\mu)\geq \lambda$, which implies that
$T_{(Id,-\dom(\mu)')} \in \Cc_G(\CalX_\lambda)$. By applying the lemma to
$-\dom(\mu)'$ for each $\mu$ such that $\dom(\mu)\geq \lambda$ we obtain the inclusion:
\[
\Cc_G(\CalX_\lambda) \subseteq \bigoplus_{\dom(\mu) \geq \lambda} \Hn
X^{\dom(\mu)'}\Hn = \Hnaffp  X^{\lambda'} \Hnaffp
\]
which gives the equality.

\end{dem}

So far we have defined for each partition $\lambda$ a closed subset
$\CalX_\lambda$ and a two sided ideal $I_\lambda=\Cc_G(\CalX_\lambda)$
generated by $X^{\lambda'}$. We can ask if every two sided ideal comes this
way.

The first remark to make is that we can associate to every finite set of
partition $\lbar=(\lambda^{(i)})_i$ the closed subset
$\bigcup_i\CalX_{\lambda^{(i)}}$. The corresponding two-sided ideal is the
sum $I_\lbar=\sum_i I_{\lambda^{(i)}}$. It is easy to see that the bicommutant
theorem holds in this case too (the proofs are the same).
 
\begin{thm}
Every $G$-stable closed subset $\CalF$ of $\CalX$ such
that $\Cc_G(\CalF)$ is a two-sided ideal of
$\Cc_G(\CalX)$ is of the form $\CalX_\lbar$, for a finite
set of partition $\lbar=(\lambda^{(i)})_i$.
\end{thm}

\begin{dem}
As $\Cc$-vector spaces, we have
\begin{equation}
\label{dec2}
\Cc_G(\CalF)=\bigoplus_{\CalO_w
  \subseteq \CalF} \Cc T_w.
\end{equation}
The next lemma is a refinement of lemma \ref{lem1}.
\begin{lem}
For every $w \in \Waffn$ we have:
\[
\Hn T_w \Hn = \bigoplus_{\sigma,\sigma' \in \Wn} \Cc T_{\sigma w \sigma'}
\]
\end{lem}

We use the lemma to rewrite the sum (\ref{dec2}) as:
\begin{equation}
\label{dec3}
\Cc_G(\CalF)= \sum_{\CalO_{(Id,\lambda)}\subseteq \CalF} \Hn
T_{(Id,\lambda)} \Hn
\end{equation}
where the sum is over the dominant partitions $\lambda$. Indeed, every $w
\in \Waffn$ belongs to a class $\Wn (Id,\lambda) \Wn$ for a dominant
$\lambda$.

We have the usual partial order on the partitions $\lambda$, and the set of
minimal partitions $(\lambda^{(i)})=\lbar$ is finite. Using that
$\Cc_G(\CalF)$ is in fact a ($\Hnaffp$,$\Hnaffp$)-bimodule, the equality (\ref{dec3})
gives:
\[
\Cc_G(\CalF)= \sum_i \Hnaffp T_{(Id,\lambda^{(i)})} \Hnaffp = \Cc_G(\CalX_\lbar)
\]
Then $\CalF=\CalX_\lbar$.
\end{dem}

\begin{rem}
There is an established Schur-Weyl duality between cyclotomic Hecke
algebras and the so-called cyclotomic $q$-Schur algebras (cf \cite{SS},
\cite{A}), but only in the semi-simple case. Our quotients are cyclotomic
Hecke algebras (with all parameters equal to zero) when we take the
partition \linebreak $(d,0,\cdots,0)$ and the semi-simplicity is obviously not
verified in this case.
\end{rem}

\section{Canonical basis of $\Hnaffl$}

The geometric construction of our algebras allows us to construct canonical
bases for them. Such bases, which are also called Kazhdan-Lusztig bases for
Hecke algebras, were introduced for quantum enveloping algebras by
Kashiwara and Lusztig (see \cite{L1}). These bases have several important
properties, which include positivity of the structure constants and
compatibility with bases of representations (see \cite{A2},
\cite{LLT},\cite{VV}).

Write $\zeta:\Hnaffp \rightarrow \Hnaffl$ for the quotient map.

Let us call $B$ the canonical basis of the affine Hecke algebra
$\Hnaff$. This basis $B=(b_\CalO)$ is defined by the formula:
\[
b_\CalO= \sum_{i,\CalO'} v^{-i+dim\CalO} \mathrm{dim}\mathcal{H}^i_{\CalO'}(IC_\CalO)1_{\CalO'}
\]
where $\mathcal{H}^i_{\CalO'}(IC_\CalO)$ is the fiber at any point in
$\CalO'$ of the cohomology sheaf of the intersection complex of $\CalO$.

As $\CalX$ is a closed subset of $\Xchap \times \Xchap$, the subset $B^+$
of $B$ defined by $B^+=\{ b \in B, b \in \Cc_G(\CalX) \}$ is a basis of $\Cc_G(\CalX)$.

\begin{thm}
The set of elements:
\[
B'=\{\zeta (b), b\in B^+ \ \arrowvert  \ \zeta(b)\neq0 \}
\]
form a basis of $\Hnaffl$.
\end{thm}

\begin{dem}
It suffices to see that:
\[
\Ker(\zeta)= \bigoplus_{\phi(b)=0} \Cc b
\]
 As $\Ker(\zeta)$ is the set of functions supported on the closed subset
$\CalX_\lambda$, the elements $b_\CalO$, where $\CalO \subseteq
\CalX_\lambda$ form a basis of $\Ker(\zeta)$. The theorem follows.

\end{dem}

\section{The case $d<n$}

The previous bicommutant theorems are true only in the case $d\geq n$. In
the case $d<n$, one half of the result still holds.

\begin{thm}
If $d<n$ the map:
\[
\Cc_G(X \times X) \rightarrow \End_{\Cc_G(Y \times Y)}(\Cc_G(Y \times X))
\]
is surjective, when $X$ is the complete (resp. affine) flag variety  and $Y$
the (resp. affine) flag variety of length $d$.
\end{thm}

\begin{dem}
Let $d<n$. We associate as before to each composition $\dbar$ of $n$ of length $d$ with a
connected component $Y_\dbar$ of the affine flag variety of length
$d$. To a composition of $n$ of length $d$ we associate a subset of $S=
\{0,\cdots,n-1\}$ of order at least $n-d$, in the way that the composition
of $n$ give the sequence dimensions of successive factors in the flag while
the set $I$ give which step in a complete flag are forgotten. We have a
bijection between these subsets and isomorphism classes of connected
components in $Y$.

We then write $Y_I$ for a connected component in the corresponding class.

Let $W$ be the extended affine Weyl group of $GL_n$. We recall that
it is the semi-direct product $W'\rtimes \Omega$, where $W'$ is the affine
Weyl group of $GL_n$ and $\Omega$ is isomorphic to $\Zz$, generated by
$\rho$. The group $W'$ is a Coxeter group, which is equipped with the
length function $l$. Each element $w$ of $W$ can be uniquely written
$w'\rho^z$, where $w'$ is an element of $W'$ and $z\in\Zz$. We define
the length $l(w)$ of $w$ by $l(w')$ and its height $h(w)=|z|$.

For each $I\subseteq S$ there is a $G$-invariant surjective map:
\[
\phi_I: X \twoheadrightarrow Y_I
\]
and for each $I\subseteq J$ there is also a surjective $G$-morphism:
\[
\phi_{I,J}:Y_I \twoheadrightarrow Y_J
\]
From these we deduce the maps:
\[
\theta_I: \Cc_G(X \times X) \rightarrow \Cc_G(Y_I \times X)
\]
\[
\theta_{I,J}: \Cc_G(Y_I \times X) \rightarrow \Cc_G(Y_J \times X)
\]

given by:
\[
\theta_I(f)(L,L')= \sum_{L'' \in \phi_I^{-1}(L)} f(L'',L') =
1_{\Delta(Y_I \times X)}*f
\]
and
\[
\theta_{I,J}(g)(L,L')= \sum_{L'' \in \phi_{I,J}^{-1}(L)}
g_I(L'',L)= 1_{\Delta(Y_J \times Y_I)}*g_I
\]

They commute with the right action of $\Cc_G(X\times X)$ by convolution.

By summing over the $I\subseteq S$ of order greater or equal to $n-d$, we
define a map:
\[
\theta: \Cc_G(X\times X)  \rightarrow \Cc_G(Y\times X)
\]

\begin{lem}\label{d<n}
The image of the map $\theta$ is the set:
\[
\{ f= \sum_{|I|\geq n-d} f_I \in \Cc_G(Y\times X) \ \arrowvert \  \forall I \subseteq J, \
f_J= \theta_{I,J}(f_I) \}
\]
\end{lem}
\begin{dem}

The inclusion of the image of $\theta$ in this set comes from the equality
$\theta_{I,J} \circ \theta_I = \theta_J$.

Let's prove the other inclusion. We must find, given a family of functions
$(f_I)_{|I|\geq n-d}$ verifying for each $I\subseteq J$ the equality
$\theta_{I,J}(f_I)=f_J$, a function $f$ in $\Cc_G(X\times X)$ such that
$f_I=\theta_I(f)$ for each $I$.

In order to give a function in $\Cc_G(X \times X)$, we have to give a value
for each orbit $\CalO_w$, with $w \in W$. By abuse of notation, we will
denote that value by $f(\CalO_w)$. We proceed in two steps.

Consider the set $M$ of all $w \in W$ such that for each $I$ of order at least
$n-d$, the orbit $\CalO_w$ is not open in the fiber $(\phi_I
\times Id)^{-1}(\CalO_{W_Iw})$, where $\CalO_{W_Iw}$ is the image of $\CalO_w$
in $Y_I \times X$. It is equivalent to say that for each $I$ of order at
least $n-d$ (and strictly less than $n$), the element $w$ is not of maximal
length in the class $W_Iw$ (seen as a subset of $W$),
where $W_I$ is the Young subgroup of $W$ generated by the elements $s_i$
with $i \in I$ (it is a finite group because $|I| < n$).

The functions $f_I$ have a compact support, so we can choose $k$ such
that for each $w\in W$ of length or height greater or equal to $k$, $f_I$
vanishes on $\CalO_{W_I w}$ for every $I$.

Fix an integer $l>k$.
For each element $w$ of $M$ of length less or equal to $l$ we assign an
arbitrary value to $f(\CalO_w)$, and for the element of length or height  greater or equal
to $l$ we set $f(\CalO_w)=0$.

Now we have to give a value to to the elements $w \in W - M$. By
definition, for such an element $w$ there exists a set $I$ for which the orbit
$\CalO_w$ is dense in the fiber $(\phi_I \times
Id)^{-1}(\CalO_{W_Iw})$.

We know that $w$ is of maximal length in $W_Iw$ if and only if for each $i
\in I$ we have $l(s_iw)<l(w)$. We can deduce from this that there is a
maximal set $I(w)$ such that $w$ is of maximal length in $W_{I(w)}w$. 

Now we give a value to the elements $w \in M$, proceeding by induction on
the length of $w$. To that purpose we use the following equation, where we
write simply $I$ for $I(w)$:
\[
m_I(w)f(\CalO_w)=f_I(\CalO_{W_Iw})-\sum_{\substack{w' \in
    W_Iw \\ w' \neq w}} m_I(w')f(\CalO_{w'})
\]
where $m_I(w')= \arrowvert \{ L_\bullet'' \in \phi_I^{-1}(L_\bullet), \
(L_\bullet'',L_\bullet') \in \CalO_{w'} \} \arrowvert$ for any $L_\bullet
\in Y_I$ such that $(L_\bullet,L_\bullet') \in \CalO_{W_Iw}$.

This determines the values $f(\CalO_w)$ because each element in
the sum has a length strictly smaller than $l(w)$, so that $f(\CalO_{w'})$ is
already defined.

It remains to show that beyond a fixed length or height the obtained values
$f(\CalO_w)$ are zero (for the function to be compactly supported), and
that the function given by this method is a solution to our problem.

We start with the first point.

We prove first that for each element $w$ of length greater or equal to
$l+\frac{n(n-1)}{2}|I(w)|$, $f(\CalO_w)$ is zero.

 If $w$ is not maximal in any of its classes $W_Iw$, then we have taken
 $f(\CalO_w)=0$. If not, for $I=I(w)$, the element $w$ is maximal in its
 class $W_Iw$, and $f(\CalO_w)$ is given by:
\[
m_I(w)f(\CalO_w)=f_I(\CalO_{W_Iw})-\sum_{\substack{w' \in
    W_Iw \\ w' \neq w}} m_I(w')f(\CalO_{w'})
\]

In the sum the elements $w'$ have a length greater or equal to
$l+\frac{n(n-1)}{2}|I(w)|-\frac{n(n-1)}{2}$. But as they are not maximal in $W_Iw$, if they are in $M$
they satisfy $|I(w')|<|I(w)|$. Then if $w'$ is in $M$ it has a length greater
or equal to $l+\frac{n(n-1)}{2}|I(w')|$, and if $w'$ is not in $M$ we have
$f(\CalO_{w'})=0$. By induction on $|I(w')|$, we easily see that each $f(\CalO_{w'})$ in
the right sum is zero, and $f_I(\CalO_{W_Iw})$ too, so $f(\CalO_w)$ is
zero. We also know that $f(\CalO_w)$ is zero for each element $w$ of length
greater or equal to $l$. Therefore $f$ is non-zero only on a finite set.

It remains to check that the function just built is a solution to our
problem. We have to show that for each $I$ and for each $w\in W$, the
following holds:
\[
f_I(\CalO_{W_Iw})=\sum_{w'\in W_Iw} m_I(w')f(\CalO_{w'})
\]
Obviously, it is sufficient to check this equation when $w$ is of maximal
length in its class $W_Iw$. If $I=I(w)=:J$, then the equation is true by
construction of $f(\CalO_w)$. If not we have $I \subseteq J$. We proceed by
induction: suppose that the equality is true for each $x$ such that
$l(x)<l(w)$.

By construction of $f(\CalO_w)$ we have:
\[
f_J(\CalO_{W_Jw})=\sum_{w''\in W_Jw} m_J(w'')f(\CalO_{w''})
\]
If we decompose the sum along the elements of the class $W_I \setminus W_J
w$ we obtain:
\[
f_J(\CalO_{W_Jw})=\sum_{x\in W_I \setminus W_Jw} \sum_{z\in W_Ix} m_J(z)f(\CalO_z)
\]
\[
f_J(\CalO_{W_Jw})=\sum_{\substack{x\in W_I \setminus W_Jw \\ x
    \neq W_Iw}} \sum_{z\in W_Ix}
    m_J(z)f(\CalO_z)+\sum_{y\in W_Iw} m_J(y)f(\CalO_y)
\]
For $w' \in W$ we define $m_{I,J}(W_Iw')= \arrowvert \{ L_\bullet'' \in
\phi_{I,J}^{-1}(L_\bullet), \  (L_\bullet'',L_\bullet')\in \CalO_{W_Iw'} \}
\arrowvert$ for any $L_\bullet \in Y_J$ such that $(L_\bullet,L_\bullet')
\in \CalO_{W_Jw}$. From the identity $\phi_J=\phi_{I,J} \phi_I$ we deduce
$m_J(w')=m_{I,J}(W_Iw') m_I(w')$. So we can write:
\[
f_J(\CalO_{W_Jw})=\sum_{\substack{x\in W_I \setminus W_Jw \\ x
    \neq W_Iw}} m_{I,J}(W_Ix)\sum_{z\in W_Ix}
    m_J(z)f(\CalO_z)+ m_{I,J}(W_Iw)\sum_{y\in W_Iw} m_I(y)f(\CalO_y)
\]

(where $W_I \setminus W_Jw$ is seen as a subset of the quotient $W_I
\setminus W$)

 But each $z \in W_Ix$ for $x\neq W_Iw$ is of length less than
 $l(w)$. Then we can use 
\[
\sum_{z\in W_Ix} m_I(z)f(\CalO_z)=f_I(\CalO_{W_Ix})
\]
So we have:
\[
f_J(\CalO_{W_Jw})=\sum_{\substack{x\in W_I \setminus W_Jw \\ x
    \neq W_Iw}} m_{I,J}(W_Ix)f_I(\CalO_{W_Ix})+\sum_{y\in W_Iw} m_I(y)f(\CalO_y)
\]
But from the equality $f_J=\theta_{I,J}(f_I)$ the following holds:
\[
f_J(\CalO_{W_Jw})=\sum_{\substack{x\in W_I \setminus
    W_Jw \\ x \neq W_Iw}} m_{I,J}(W_Ix) f_I(\CalO_{W_Ix})+ m_{I,J}(W_Iw)f_I(\CalO_{W_Iw})
\]
We finally obtain:
\[
f_I(\CalO_{W_Iw})=\sum_{x\in W_I w} m_I(x)f(\CalO_x)
\]

\end{dem}

Let $P$ be an element of $\End_{\Cc_G(Y\times Y)}(\Cc_G(Y\times X))$. As $P$
commmutes with the action of the characteristic functions of the diagonals
of the components $Y_I\times X$, the subspaces  $\Cc_G(Y_I \times X)$ are
stable for the endomorphism $P$. So we can write $P= \oplus_I P_I$, where
$P_I \in \End_{\Cc_G(Y_I \times Y_I)}(Y_I \times X)$.

Write $\psi_I:\Cc_G(Y_I\times Y_I) \hookrightarrow \Cc_G(Y_I \times X)$
for the injection deduced from the surjection $\phi_I:X \rightarrow
Y_I$. Write $f_I:=P_I(1_{\Delta(Y_I \times X)})$. For each
$g\in \Cc_G(Y_I \times X)$ the following equalities hold:
\[
P_I(g)=P_I(\psi_I(g)* 1_{\Delta(Y_I\times
  X)})=\psi_I(g)*P_I(1_{\Delta(Y_I\times X)})=\psi_I(g)*f_I
\]
The next step is to apply the lemma \ref{d<n} to lift the $f_I$s to some
$f\in \Cc_G(X \times X)$.

For this we have to check that the functions $f_I=P_I(1_{\Delta(Y_I \times X)})$
satisfy the hypothesis of the lemma.

But for $I \subseteq J$, we have $1_{\Delta(Y_J \times Y_I)}*1_{\Delta(Y_I \times
  X)}=1_{\Delta(Y_J \times X)}$, thus $f_J= \theta_{I,J}(f_I)$.

We deduce that there exists $f \in \Cc_G(X \times X)$ such that for
each $I$ of order greater or equal to $n-d$ we have $f_I=\theta_I(f)$. But
by construction for $g \in \Cc_G(Y_I \times X)$ the product $\psi_I(g)*f_I$
is equal to $g*f$. We have shown that $P(g)=g*f \ , \ \forall g \in \Cc_G(Y
\times X)$. The theorem follows.

\end{dem}

\end{document}